\documentstyle[12pt]{article}

\begin{document}
\newcommand{\ol }{\overline}
\newcommand{\ul }{\underline }
\newcommand{\ra }{\rightarrow }
\newcommand{\lra }{\longrightarrow }
\newcommand{\ga }{\gamma }
\newcommand{\st }{\stackrel }
\newcommand{\scr }{\scriptsize }
\title{\Large\textbf{Some Baer Invariants of Free Nilpotent Groups}}
\author{\textbf{Behrooz Mashayekhy} \\
Department of Mathematics, Ferdowsi University of Mashhad,\\
P.O.Box 1159-91775, Mashhad, Iran, and \\ Institute for Studies in
Theoretical Physics and Mathematics (IPM)\\ \textbf{Mohsen
Parvizi}\\ Department of Pure Mathematics, Damghan University of
Basic
Sciences,\\ P.O.Box 36715-364, Damghan, Iran.\\
E-mail: mashaf@math.um.ac.ir   \&  parvizi@dubs.ac.ir}
\date{ }
\maketitle
\begin{abstract}
We present an explicit structure for the Baer invariant of a free
$n$th nilpotent group (the $n$th nilpotent product of infinite
cyclic groups, $\textbf{ Z}\st{n}* \textbf{ Z}\st{n}*\ldots
\st{n}*\textbf{ Z}$) with respect to the variety ${\cal V}$ with the
set of words $V=\{[\ga_{c_1+1},\ga_{c_2+1}]\}$, for all $c_1\geq
c_2$ and $2c_2-c_1>2n-2$. Also, an explicit formula for the
polynilpotent multiplier of a free $n$th nilpotent group is given
for any class row $(c_1,c_2,\ldots,c_t)$, where $c_1\geq n$.
\end{abstract}
\newpage
\begin{center}
\textbf{1. Introduction and Preliminaries}\\
\end{center}

 One of the interesting problem, related to the well-known notion
 of the Schur multiplier and its varietal generalization, the Baer invariant, is to
 obtain some structures  for the Baer invariant of a given group $G$ and specially of
 some famous products of groups, such as direct products, semidirect products, free products,
 nilpotent products and regular products. Determining
these Baer invariants of a given group is known to be very useful
for classification of groups into isologism classes. Also,
structures of Baer invariants are very essential for studying
generalized capability and covering groups. Some people have
succeeded to find the structures
 as follows.

I. Schur [18] in 1907 and J. Wiegold [20] in 1971 obtained the
structure of the Schur multiplier of the direct product of two
finite groups as follows:
\begin{center}
$M(A\times B)\cong M(A)\oplus M(B)\oplus \frac{[A,B]}{[A,B,A*B]}$,
where $\frac{[A,B]}{[A,B,A*B]}\cong A_{ab}\otimes B_{ab}$.
\end{center}
In 1979, M. R. R. Moghaddam [15] and in 1998, G. Ellis [2],
succeeded to extend the above result to obtain the structure of the
$c$-nilpotent multiplier of the direct product of two groups, ${\cal
N}_cM(A\times B)$. Also in 1997 the first author in a joint paper
[10] presented an explicit formula for the $c$-nilpotent multiplier
of a finite abelian group.

K. I. Tahara [19] in 1972 and W. Haebich [6] in 1977 found some
structures for the Schur multiplier of the semidirect product of two
groups. Also, the first author [9,13] extended some of the above
results to the variety of nilpotent groups.

In 1972 W. Haebich [5] presented a formula for the Schur multiplier
of a regular product of a family of groups. It is known that the
regular product is a generalization of the nilpotent product and the
last one is a generalization of the direct product, so Haebich's
result is an interesting generalization of the Schur's result. Also,
M. R. R. Moghaddam [16], in 1979 gave a formula similar to Haebich's
formula for the Schur multiplier of a nilpotent product. Moreover,
in 1992, N. D. Gupta and M. R. R. Moghaddam [4] presented an
explicit formula for the $c$-nilpotent multiplier of the $n$th
nilpotent product $\textbf{Z}_2\st{n}*\textbf{Z}_2$. G. Ellis [3]
remarked that there is a slip in the statement of the result in case
$c\leq n-1$ and gave the correct structure.

In 2001, the first author [11] found a structure similar to
Haebich's type for the $c$-nilpotent multiplier of a nilpotent
product of a family of cyclic groups. The $c$-nilpotent multiplier
of a free product of some cyclic groups was studied by the first
author [12] in 2002.

Recently, the authors [14,17] concentrated on the Baer invariant
with respect to the variety of polynilpotent groups, for the first
time. We presented an explicit structure for some polynilpotent
multipliers of the $n$th nilpotent product of some infinite cyclic
groups [17] and also found explicit structures for all polynilpotent
multipliers of finitely generated abelian groups [14].

Now trying to extend the above results to the vast variety of
\textit{outer commutators} we concentrate, as a first step, on a
variety $\cal V$ with a set of words $\{
[\ga_{c_1+1},\ga_{c_2+1}]\}$. In this paper we intend to extend the
last above results in two directions. First, to obtain an explicit
formula for some Baer invariants of the $n$th nilpotent product of
some infinite cyclic groups,
$${\cal V}M(\textbf{ Z}\st{n}* \textbf{ Z}\st{n}*\ldots \st{n}*\textbf{ Z}),$$
in which $\cal V$ is the above variety for all $c_1\geq c_2$ and
$2c_2-c_1> 2n-2$. Note that the first restriction on the parameters,
$c_1\geq c_2$, is very natural and the second one help us to
calculate the Baer invariant. Second, to present an explicit formula
for the polynilpotent multiplier of the $n$th nilpotent product of
some infinite cyclic groups
$${\cal N}_{c_1,\ldots ,c_t}M(\textbf{ Z}\st{n}* \textbf{ Z}\st{n}*\ldots \st{n}*\textbf{ Z}),$$
for any class row $(c_1,c_2,\ldots,c_t)$, where $n\leq c_1$.

In the following you can find some preliminaries which are used in
our method.\\
\textbf{Definition 1.1}. Let $G$ be any group with a free
presentation $G\cong F/R$, where $F$ is a free group. Then, after R.
Baer [1], the \textit{Baer invariant} of $G$ with respect to a
variety of groups ${\cal V}$, denoted by ${\cal V}M(G)$, is defined
to be
$${\cal V}M(G)=\frac{R\cap V(F)}{[RV^*F]}\ ,$$
where $V$ is the set of words of the variety ${\cal V}$, $V(F)$ is
the verbal subgroup of $F$ with respect to ${\cal V}$ and
$$[RV^*F]=<v(f_1,\ldots ,f_{i-1},f_ir,f_{i+1},\ldots,f_n)v(f_1,\ldots,f_i,
\ldots,f_n)^{-1}\mid $$
$$r\in R, 1\leq i\leq n, v\in V ,f_i\in F, n\in {\bf N}>.$$

In special case of the variety ${\cal A}$ of abelian groups, the
Baer invariant of $G$ will be the well-known notion the
\textit{Schur multiplier}
$$\frac{R\cap F'}{[R,F]}.$$

If ${\cal V}$ is the variety of nilpotent groups of class at most
$c\geq1$, ${\cal N}_c$, then the Baer invariant of $G$ with
respect to ${\cal N}_c$ which is called the \textit{$c$-nilpotent
multiplier} of $G$, will be
$${\cal N}_cM(G)=\frac{R\cap \gamma_{c+1}(F)}{[R,\ _cF]},$$
where $\gamma_{c+1}(F)$ is the $(c+1)$-st term of the lower
central series of $F$ and $[R,\ _1F]=[R,F], [R,\ _cF]=[[R,\
_{c-1}F],F]$, inductively.\\ \ \ \\
\textbf{Lemma 1.2} (J. A. Hulse and J. C. Lennox 1976). \textit{If
$u$ and $w$ are any two words and $v=[u,w]$ and $K$ is a normal
subgroup of a group $G$, then}
$$ [Kv^*G]=[[Ku^*G],w(G)][u(G),[Kw^*G]].$$
\textit{Proof.} See [8, Lemma 2.9].

 Now, using the above lemma, let $\cal V$ be the outer commutator variety of
 groups defined by the set of words $\{ [\ga_{c_1+1},\ga_{c_2+1}]\}$, then
 the Baer invariant of a group $G$
 with respect to this variety, is as follows:
 $${\cal V}M(G)\cong \frac{R\cap
 [\ga_{c_1+1}(F),\ga_{c_2+1}(F)]}{[R,\ _{c_1}F,\ga_{c_2+1}(F)][R,\ _{c_2}F,\ga_{c_1+1}(F)]},
 \ \ (\star).$$
\textbf{Definition 1.3}. \textit{Basic commutators} are defined in
the usual way. If $X$ is a fully ordered independent subset of a
free group, the basic
commutators on $X$ are defined inductively over their weight as follows:\\
$(i)$ the members of $X$ are basic commutators of weight one on $X$;\\
$(ii)$ assuming that $n>1$ and that the basic commutators of
weight less than $n$ on $X$ have been defined and ordered;\\
$(iii)$ a commutator $[b,a]$ is a basic commutator of weight $n$ on
$X$ if  $wt(a)+wt(b)=n,\ a<b$, and if $b=[b_1,b_2]$, then $b_2\leq
a$. The ordering of basic commutators is then extended to include
those of weight $n$ in any way such that those of weight less than
$n$ precede those of weight $n$. The natural way to define the order
on basic commutators of the same weight is lexicographically,
$[b_1,a_1]<[b_2,a_2]$ if $b_1<b_2$ or if $b_1=b_2$ and $a_1<a_2$.

The next two theorems are vital in our investigation.\\
\textbf{Theorem 1.4} (P. Hall [7]). \textit{Let $F=<x_1,x_2,\ldots
,x_d>$ be a free group, then
$$ \frac {\ga_n(F)}{\ga_{n+i}(F)} \ \ , \ \ \ \  1\leq i\leq n$$
is the free abelian group freely generated by the basic
commutators of weights
$n,n+1,\ldots ,n+i-1$ on the letters $\{x_1,\ldots ,x_d\}.$}\\
\textbf{Theorem 1.5} (Witt Formula [7]). \textit{The number of basic
commutators of weight $n$ on $d$ generators is given by the
following formula:
$$ \chi_n(d)=\frac {1}{n} \sum_{m|n}^{} \mu (m)d^{n/m},$$
where $\mu (m)$ is the M\"{o}bius function.} \\
\textbf{Definition 1.6}. Let $\cal V$ be a variety of groups defined
by a set of laws $V$. Then the \textit{verbal product} of a family
of groups $\{G_i\}_{i\in I}$ associated with the variety $\cal V$ is
defined to be
$${\cal V}\prod_{i\in I}G_i=\frac {\prod^*G_i}{V(G)\cap
[G_i]^*},$$ where $G=\prod^{*}_{i\in I}G_i$ is the free product
of the family $\{G_i\}_{i\in I}$ and $[G_i]^*=\\ <[G_i,G_j]|i,j\in
I,i\neq j>^G$ is the cartesian subgroup of the free product $G$.

 The verbal product is also known as \textit {varietal product} or
 simply \textit{$\cal V$-product}. If $\cal V$ is the variety of
 all groups, then the corresponding verbal product is the free
 product; if ${\cal V}={\cal A}$ is the variety of all abelian
 groups, then the verbal product is the direct product. Also, if
 ${\cal V}={\cal N}_c$ is the variety of nilpotent groups of class
 at most $c\geq 1$, then the verbal product is called the
 $c$\textit{th nilpotent product} of the $G_i$'s.\\

\begin{center}
\textbf{2. The Main Results}\\
\end{center}
Let $G\cong \textbf{ Z}\st{n}* \textbf{ Z}\st{n}*\ldots
\st{n}*\textbf{ Z}$ be the $n$th nilpotent product of $m$ copies of
the infinite cyclic group $\textbf{ Z}$. Using Definition 1.6, it is
easy to see that $G$ is the free $n$th nilpotent group of rank $m$,
and so has the following free presentation
$$1\lra \ga_{n+1}(F)\lra F\lra G\lra 1,$$
where $F$ is the free group on a set $X=\{x_1,x_2,\ldots,x_m\}$.

Now, we try to obtain the structure of some outer commutator
multipliers of $G$ of the form
$${\cal V}M(G),$$
where ${\cal V}$ is defined by the set of words $\{
[\ga_{c_1+1},\ga_{c_2+1}]\}$, $c_1\geq c_2$, and $2c_2-c_1>2n-2$
(Note that $\ga_c=[x_1,\ldots, x_c]$).

Using $(\star)$ we have $${\cal V}M(G)\cong \frac{\ga_{n+1}(F)\cap
 [\ga_{c_1+1}(F),\ga_{c_2+1}(F)]}{[\ga_{c_1+n+1}(F),\ga_{c_2+1}(F)][\ga_{c_2+n+1}(F),
 \ga_{c_1+1}(F)]}.$$
Now if we have $c_1+c_2+1\geq n$, then

$${\cal V}M(G)\cong \frac{[\ga_{c_1+1}(F),\ga_{c_2+1}(F)]}{[\ga_{c_1+n+1}(F),\ga_{c_2+1}(F)]
[\ga_{c_2+n+1}(F),\ga_{c_1+1}(F)]}.$$

 In order to find the structure of ${\cal V}M(G)$, we
need the following notation and lemmas. Using Definition and
Notation 1.2, we defined the following sets when $c_1\geq c_2$.
\begin{center}
$A=\{[\beta,\alpha] \ | \ \beta $ and $ \alpha $ are basic
commutators on $X$ such that $ \beta>\alpha, $ $ \ c_1+1\leq
wt(\beta)\leq c_1+n, \   c_2+1\leq wt(\alpha)\leq c_2+n \};$\\ \ \  \\

$B=\{ \ [\beta,\alpha] \ | \ \beta$ and $\alpha$ are basic
commutators on $X$ such that $\beta>\alpha$, $c_1+n+1\leq
wt(\beta)$ , $c_2+1\leq wt(\alpha)$ , $wt(\beta)+wt(\alpha)\leq
2n+c_1+c_2+1 \};$\\ \ \ \\

$C=\{ \ [\beta,\alpha] \ | \ \beta$ and $\alpha$ are basic
commutators on $X$ such that $\beta>\alpha$, $c_2+n+1\leq
wt(\beta)$ , $c_1+1\leq wt(\alpha)$ , $wt(\beta)+wt(\alpha)\leq
2n+c_1+c_2+1 \}$.
\end{center}
\textbf{ Lemma 2.1}. \textit{If $n-1\leq 2c_2-c_1$,
then every element of $A$ is a basic commutator on $X$}.\\
\textit{ Proof.} Every element of $A$ has the form
$[\beta,\alpha]$, where $\beta$ and $\alpha$ are basic commutators
on $X$, $\beta>\alpha$ and $c_1+1\leq wt(\beta)\leq c_1+n, \
c_2+1\leq wt(\alpha)\leq c_2+n$. Now, let
$\beta=[\beta_1,\beta_2]$, then in order to show that
$[\beta,\alpha]$ is a basic commutator on $X$, it is enough to
show that $\beta_2\leq \alpha$. Since $\beta=[\beta_1,\beta_2]$ is
a basic commutator on $X$, so $\beta_1>\beta_2 $ and hence
$wt(\beta_2)\leq \frac{1}{2}wt(\beta)$. Now, if $n-1\leq
2c_2-c_1$, then $\frac{1}{2}(c_1+n)<c_2+1$. Thus, since
$wt(\beta)\leq c_1+n$, we have
$$wt(\beta_2)\leq \frac{1}{2}wt(\beta)\leq \frac{1}{2}(c_1+n)<c_2+1\leq wt(\alpha).$$
Therefore $\beta_2<\alpha$ and hence the result holds. $\Box$\\
\textbf{ Lemma 2.2}. \textit{With the above notation,\\
$(i)$ if $2n-2<2c_2-c_1$, then every
element of $B$ is a basic commutator on $X$;\\
$(ii)$ if $2n-2<2c_1-c_2$, then every
element of $C$ is a basic commutator on $X$.}\\
\textit{ Proof.}$(i)$ Let $[\beta,\alpha]$ be an element of $B$.
By definition of $B$, $\beta$ and $\alpha$ are basic commutators
on $X$ and $\beta>\alpha$. Now, let $\beta=[\beta_1,\beta_2]$,
where $\beta_1$ and $\beta_2$ are basic commutators on $X$. It is
enough to show that $\beta_2\leq \alpha$. Since
$\beta=[\beta_1,\beta_2]$ is a basic commutator, so
$\beta_1>\beta_2$ and hence $wt(\beta_2)\leq
\frac{1}{2}wt(\beta)$. By definition of $B$, we have
$wt(\beta_2)\leq \frac{1}{2}wt(\beta)\leq
\frac{1}{2}(2n+c_1+c_2+1-c_2-1)=\frac{1}{2}(2n+c_1)<c_2+1\leq
wt(\alpha).$\\
$(ii)$ A similar argument can be applied for elements of $C$.$\Box$\\
\textbf{ Lemma 2.3}. \textit{With the above notation, the following statements hold.\\
$(i)$ If $c_2+n<c_1+1$, then
$$ |A|=(\sum_{i=c_1+1}^{c_1+n}\chi_i(m))(\sum_{c_2+1}^{c_2+n}\chi_i(m)).$$\\
$(ii)$ If $c_2+n\geq c_1+1$, then}
$$
|A|=(\sum_{i=c_1+1}^{c_1+n}\chi_i(m))(\sum_{c_2+1}^{c_1}\chi_i(m))+$$
$$(\sum_{i=c_2+n+1}^{c_1+n}\chi_i(m))(\sum_{c_1+1}^{c_2+n}\chi_i(m))+
\chi_2(\sum_{i=c_1+1}^{c_2+n}\chi_i(m)).$$\\
\textit{ Proof.} $(i)$ We have $c_2+1\leq wt(\alpha)\leq
c_2+n<c_1+1\leq wt(\beta)\leq c_1+n$, so $\beta>\alpha$ and hence
all possible of $\beta$ and $\alpha$ is accepted. Clearly the
number of $\beta$'s is $(\sum_{i=c_1+1}^{c_1+n}\chi_i(m))$ and
the number of $\alpha$'s
is $(\sum_{c_2+1}^{c_2+n}\chi_i(m))$. Thus the result holds.\\
$(ii)$ Let $c_1+1\leq c_2+n $. Then $c_2+1\leq c_1+1\leq
c_2+n\leq c_1+n$. Put
\begin{center}
$A_1=\{[\beta,\alpha] \ | \ \beta $ and $ \alpha $ are basic
commutators on $X$ such that $ \beta>\alpha, $ $ \ c_1+1\leq
wt(\beta)\leq c_1+n$ , $c_2+1\leq wt(\alpha)< c_1+1  \};$\\ \ \  \\

$A_2=\{ \ [\beta,\alpha] \ | \ \beta$ and $\alpha$ are basic
commutators on $X$ such that $\beta>\alpha$, $c_1+1\leq
wt(\beta)\leq c_2+n$ , $c_1+1\leq wt(\alpha)\leq c_2+n  \};$\\ \ \ \\

$A_3=\{ \ [\beta,\alpha] \ | \ \beta$ and $\alpha$ are basic
commutators on $X$ such that $\beta>\alpha$, $c_2+n< wt(\beta)\leq
c_1+n $ , $c_1+1\leq wt(\alpha)\leq c_2+n \}$.
\end{center}
Clearly $A_1, A_2, A_3$ are mutually disjoint and $ A=A_1\cup
A_2\cup A_3$. Therefore $|A|=|A_1|+|A_2|+|A_3|$. Also, it is easy
to see that
$$|A_1|=(\sum_{i=c_1+1}^{c_1+n}\chi_i(m))(\sum_{c_2+1}^{c_1}\chi_i(m));$$
$$|A_2|=\chi_2(\sum_{i=c_1+1}^{c_2+n}\chi_i(m));$$
$$|A_3|=(\sum_{i=c_2+n+1}^{c_1+n}\chi_i(m))(\sum_{c_1+1}^{c_2+n}\chi_i(m)).$$
Hence the result holds. $\Box$ \\
\textbf{ Lemma 2.4}. \textit{With the above notation and
assumption,\\
$(i)$ if $c_2+n< c_1+1$, then $A\cap C=\emptyset$;\\
$(ii)$ if $c_2+n\geq c_1+1$ , then $|A\cap
C|=(\sum_{i=c_2+n+1}^{c_1+n}\chi_i(m))(\sum_{c_1+1}^{c_2+n}\chi_i(m)).$}\\
\textit{ Proof.} By definition of $A$ and $C$ we have\\
$ A\cap C=\{\ [\beta,\alpha] \ | \ \beta$ and $\alpha$ are basic
commutators on $X$ such that $\beta>\alpha$, $\max\{c_1+1,
c_2+n+1\}\leq wt(\beta)\leq c_1+n $, $\max\{c_1+1, c_2+1\}\leq
wt(\alpha)\leq c_2+n \}.$ \\
$(i)$ Since $c_2\leq c_1$, $c_1+1\leq wt(\alpha)\leq c_2+n$ which
is a contradiction to the assumption $c_2+n<c_1+1$. Hence in this
case we have $A\cap C=\emptyset $.\\
$(ii)$ If $c_2+n\geq c_1+1$ , then $c_1+1\leq wt(\alpha)<
c_2+n+1\leq  wt(\beta)\leq c_1+n$. Thus we have always
$\beta>\alpha$ and hence the result holds. $\Box$

 The following corollary is an immediate consequence of the last
 two lemmas.\\
 \textbf{Corollary 2.5}. \textit{With the above notation and
 assumption,\\
 $(i)$ if $c_2+n<c_1+1$, then
 $$|A-C|=(\sum_{i=c_1+1}^{c_1+n}\chi_i(m))(\sum_{c_2+1}^{c_2+n}\chi_i(m));$$
 $(ii)$ if $c_2+n\geq c_1+1$, then $$|A-C|=
 (\sum_{i=c_1+1}^{c_1+n}\chi_i(m))(\sum_{c_2+1}^{c_1}\chi_i(m))+
 \chi_2(\sum_{i=c_1+1}^{c_2+n}\chi_i(m)).$$}\\
\textbf{ Lemma 2.6}. \textit{With the above notation and assumption,
if $n-1\leq c_2\leq c_1$, then we have}
$$[\ga_{c_1+1}(F),\ga_{c_2+1}(F)]\equiv <A-C>\ \  mod\  H, $$
where
$H=[\ga_{c_1+n+1}(F),\ga_{c_2+1}(F)][\ga_{c_2+n+1}(F),\ga_{c_1+1}(F)].$\\
\textit{ Proof.} Let $[\alpha,\beta]$ be a generator of
$[\ga_{c_1+1}(F),\ga_{c_2+1}(F)]$, so ${\alpha} \in
\ga_{c_+1}(F)$ and $\beta \in \ga_{c_2+1}(F)$. By P. Hall's
Theorem (1.3) we can put
$\alpha=\alpha_1\alpha_2\ldots\alpha_t\eta$ \ and \
$\beta=\beta_1\beta_2\ldots\beta_s\mu$, where
$\alpha_1,\ldots,\alpha_t$  are basic commutators on $X$ of
weights $c_1+1,\ldots, \ c_1+n$, and $\beta_1,\ldots,\beta_s$ are
basic commutators of weights $c_2+1,\ldots,c_2+n$ on $X$ and
$\eta\ \in\ga_{c_1+n+1}(F)$ and $\mu \in \ga_{c_2+n+1}(F)$.\\
By using commutator calculus, it is easy to see that
$$[\alpha,\beta]=\Pi_{i,j} [\alpha_i,\beta_j]^{f_{ij}}[\alpha_i,\mu]^{g_i}[\eta,\beta_j]^{h_j}
,$$ where $f_{ij},g_i,h_j\in\ga_{n}(F)$ (Note that,  $n-1\leq
c_2\leq c_1$). Now, it is easy to see that
\begin{center}
$[\beta_i,\eta]\in [\ga_{c_1+n+1}(F),\ga_{c_2+1}(F)]$\\
$[\mu,\alpha_i]\in [\ga_{c_2+n+1}(F),\ga_{c_1+1}(F)]$\\
$[\beta_i,\alpha_j,f_{ij}]\in
[\ga_{c_1+n+1}(F),\ga_{c_2+1}(F)][\ga_{c_2+n+1}(F),\ga_{c_1+1}(F)]$
\end{center}
(Note that the last holds by the three subgroups lemma).

Therefore we have
$$[\alpha,\beta]\equiv \Pi_{i,j}[\alpha_i,\beta_j] \ \ (mod\  [\ga_{c_1+n+1}(F),\ga_{c_2+1}(F)]
[\ga_{c_2+n+1}(F),\ga_{c_1+1}(F)]).$$
 Note that if $\alpha_i<\beta_j$, then we have $c_1+1\leq
 wt(\alpha_i)\leq wt(\beta_j)\leq c_2+n$ and so we can consider
 $[\alpha_i,\beta_j]=[\beta_j,\alpha_i]^{-1}$, where $ c_1+1\leq
 wt(\beta_j)\leq c_1+n$ and $c_2+1\leq wt(\alpha_i)\leq c_2+n$.
 Thus we have $[\beta_j,\alpha_i]\in A$ and so
it is easy to see that $\Pi_{i,j}[\alpha_i,\beta_j]\in <A>$.
Since $C\subseteq [\ga_{c_2+n+1}(F),\ga_{c_1+1}(F)]\subseteq H$,
hence
the result holds.$\Box$\\
\textbf{ Lemma 2.7}. \textit{With the above notation
and assumption, if $n-1\leq c_2\leq c_1$, then } \\
$(i) \ \ [\ga_{c_1+n+1}(F),\ga_{c_2+1}(F)]\equiv <B> \ mod \
\ga_{c_1+c_2+2n+2}(F);$\\
$(ii)  \ [\ga_{c_2+n+1}(F),\ga_{c_1+1}(F)]\equiv <C> \ mod \
\ga_{c_1+c_2+2n+2}(F)$.\\
\textit{ Proof.} $(i)$ Let $[\beta,\alpha]$ be a generator of
$[\ga_{c_1+n+1}(F),\ga_{c_2+1}(F)]$, so, $\beta\in \ga_{c_1+n+1}(F)$
and $\alpha\in \ga_{c_2+1}(F)$. By P. Hall's
Theorem (1.3) we can write\\
$\alpha=\alpha_1\alpha_2\ldots\alpha_t\eta$ and
$\beta=\beta_1\beta_2\ldots\beta_s\mu$, where
$\alpha_1,\ldots,\alpha_t$ are basic commutators of weights
$c_2+1,\ldots,c_2+n$ on $X$ and $\eta\in\ga_{c_2+n+1}(F)$, and
$\beta_1,\ldots,\beta_s$ are basic commutators of weights
$c_1+n+1,\ldots,c_1+2n$ on $X$ and $\mu\in\ga_{c_1+2n+1}(F)$.\\
By using commutator calculus, it is easy to see that
$$[\beta,\alpha]=\Pi_{i,j}[\beta_i,\alpha_j]^{f_{ij}}[\mu,\alpha_j]^{g_i}[\beta_i,\eta]^{h_j},$$
where $f_{ij},\ g_i,\ h_j\in \ga_{c_2+1}(F).$ Now we have
$$wt(\alpha_j)+wt(\mu)\geq(c_2+1)+(c_1+2n+1)\geq c_1+c_2+2n+2 $$
$$wt(\eta)+wt(\beta_j)\geq (c_2+n+1)+(c_1+n+1)\geq c_1+c_2+2n+2 $$
\begin{center}
$wt(\alpha_j)+wt(\beta_i)+wt(f_{ij})\geq
(c_2+1)+(c_1+n+1)+(c_2+1)\geq c_1+c_2+2n+2$  for all $1\neq
f_{ij}\in \ga_{c_2+1}(F)$ since $c_2\geq n-1$.
\end{center}
Therefore, by a similar argument at the end of the proof of the
previous lemma, we have
$$[\beta,\alpha]=\Pi_{i,j}[\beta_i,\alpha_j]\ \ mod \ \
\ga_{c_1+c_2+2n+2}(F) $$
$$ =\Pi_{wt(\alpha_j)+wt(\beta_i)\leq c_1+c_2+2n+1}[\beta_i,\alpha_j]\in B\ \
 mod \ \ \ga_{c_1+c_2+2n+2}(F). $$
$(ii)$ Similar to part $(i)$.$\Box$\\

Now, we are in a position to state and proof the first main result
of the paper.\\
\textbf{ Theorem 2.8}. \textit{With the above notation and
assumption, if $2c_2-c_1> 2n-2$ and $c_1\geq c_2$, then ${\cal
V}M(\textbf{ Z}\st{n}* \textbf{ Z}\st{n}*\ldots \st{n}*\textbf{ Z})$
is a free abelian group with the following basis:
$$D=\{ \ a[\ga_{c_1+n+1}(F),\ga_{c_2+1}(F)][\ga_{c_2+n+1}(F),\ga_{c_1+1}(F)]\ \ | \ \
a\in A-C \ \}.$$} \textit{Proof.} Note that if $c_1\geq c_2$ and
$2c_2-c_1>2n-2$, then it is easy to see that the following
 inequalities hold:
$$ c_1+c_2+1\geq n,\ \ c_1\geq c_2\geq n-1,\ \
2c_2-c_1\geq n-1,\ \ and\ \ 2c_1-c_2>2n-2.$$
 Therefore all of the previous lemmas hold. Clearly $${\cal
V}M(\underbrace{\textbf{ Z}\st{n}* \textbf{ Z}\st{n}*\ldots
\st{n}*\textbf{
Z}}_{m-copies})=\frac{[\ga_{c_1+1}(F),\ga_{c_2+1}(F)]}{[\ga_{c_1+n+1}(F),\ga_{c_2+1}(F)]
[\ga_{c_2+n+1}(F),\ga_{c_1+1}(F)]}$$ is an abelian group which is
generated by $D$, using Lemma 2.6. So it is enough to show that
$D$ is linearly independent. Consider the free abelian group
$\ga_{c_1+c_2+2}(F)/\ga_{c_1+c_2+2n+2}(F)$, with the basis of all
basic commutators of weights $c_1+c_2+2,\ldots, c_1+c_2+2n+1$ and
the following fact
$$\frac{H\ga_{c_1+c_2+2n+2}(F)}{\ga_{c_1+c_2+2n+2}(F)}=<\hat w\ | \ w\in B\cup C>\leq
\frac{\ga_{c_1+c_2+2}(F)}{\ga_{c_1+c_2+2n+2}(F)}.$$ Now, suppose
$$\sum_{i=1}^{\ell}k_i\bar a_i=0\ in\
\frac{[\ga_{c_1+1}(F),\ga_{c_2+1}(F)]}{[\ga_{c_1+n+1}(F),\ga_{c_2+1}(F)][\ga_{c_2+n+1}(F),
\ga_{c_1+1}(F)]},$$ where $a_1,\ldots,a_{\ell}\in A-C$, $k_i\in
\textbf{Z}.$ So we have $$\sum_{i=1}^{\ell}k_ia_i\in
[\ga_{c_1+n+1}(F),\ga_{c_2+1}(F)][\ga_{c_2+n+1}(F),\ga_{c_1+1}(F)].$$
By Lemma 2.7 the group
$[\ga_{c_1+n+1}(F),\ga_{c_2+1}(F)][\ga_{c_2+n+1}(F),\ga_{c_1+1}(F)]$
is generated by $B\cup C$ modulo $\ga_{c_1+c_2+2n+2}(F).$ Thus
considering the free abelian group
$\ga_{c_1+c_2+2}(F)/\ga_{c_1+c_2+2n+2}(F)$ we have
$\sum_{i=1}^{\ell} k_i\hat a_i=\sum_{j=1}^{t} d_j\hat w_j$ for
some $d_j\in \textbf{Z}$ and some $w_j\in B\cup C$. This implies
that $\sum k_i\hat a_i-\sum d_j\hat w_j=0$, where $a_i$'s and
$w_j$'s are basic commutators of weights $c_1+c_2+2,\ldots,
c_1+c_2+2n+1$. By the form of elements of $A$, $B$ and $C$, we
have  $(B\cup C)\cap (A-C)=\emptyset$ so $k_i=0$ and $d_j=0$, for
all $i,\ j$. Hence the result
holds.$\Box$\\
\textbf{ Corollary 2.9}. \textit{With the above notation and
assumption, if $c_1\geq c_2$ and $2c_2-c_1>2n-2$, then the
following hold.}\\
$(i)$ If $c_2+n<c_1+1$, then ${\cal V}M(\underbrace{\textbf{
Z}\st{n}* \textbf{ Z}\st{n}*\ldots \st{n}*\textbf{
Z}}_{m-copies})\cong \textbf{ Z}\oplus \ldots \oplus\textbf{ Z} \
(k-copies),$\\ where
$k=(\sum_{i=c_1+1}^{c_1+n}\chi_i(m))(\sum_{c_2+1}^{c_2+n}\chi_i(m)).$\\
$(ii)$ If $c_2+n\geq c_1+1$, then ${\cal V}M(\underbrace{\textbf{
Z}\st{n}* \textbf{ Z}\st{n}*\ldots \st{n}*\textbf{
Z}}_{m-copies})\cong \textbf{ Z}\oplus \ldots \oplus\textbf{
Z}(k'-copies), $\\
where
$k'=(\sum_{i=c_1+1}^{c_1+n}\chi_i(m))(\sum_{c_2+1}^{c_1}\chi_i(m))+
 \chi_2(\sum_{i=c_1+1}^{c_2+n}\chi_i(m)).$\\

Note that the above result is a generalization of the main result of
[17], since if $c>2n-2$, then similar formula can be obtained for
${\cal N}_{c,1}M(\textbf{ Z}\st{n}* \textbf{ Z}\st{n}*\ldots
\st{n}*\textbf{ Z})$.

 In the rest, we try to generalize the main result of [17] in
 another direction. Let $G=\textbf{ Z}\st{n}* \textbf{ Z}\st{n}*\ldots
\st{n}*\textbf{Z}$ be the free $n$th nilpotent group of rank $m$
and ${\cal N}_{c_1,\ldots ,c_t}$ be the variety of polynilpotent
groups of class row$(c_1,\ldots ,c_t)$. Consider the free
presentation $F/{\ga}_{n+1}(F)$ for $G$, where $F$ is the free
group of rank $m$. Using Lemma 1.2 we have
$${\cal N}_{c_1,\ldots ,c_t}M(G)\cong
 \frac{\ga_{n+1}(F)\cap \ga_{c_t+1}(\ldots (\ga_{c_1+1}(F))\ldots )}{[\ga_{n+1}(F),
 \ _{c_1}F,\ _{c_2}
 \ga_{c_1+1}(F),\ldots , \ _{c_t}\ga_{c_{t-1}+1}(\ldots (
 \ga_{c_1+1}(F))\ldots )]}.$$
 If $c_1\geq n$, then we can consider
 $\ga_{c_1+1}(F)/[\ga_{n+1}(F),\ _{c_1}F]$ as a free presentation
 for ${\cal N}_{c_1}M(G)$ thus we have
$${\cal N}_{c_t}M(\ldots ({\cal N}_{c_1}M(G))\ldots )\cong$$
 $$ \frac{\ga_{c_t+1}(\ldots (\ga_{c_1+1}(F))\ldots )}{[\ga_{n+1}(F),\ _{c_1}F,\ _{c_2}
 \ga_{c_1+1}(F),\ldots , \ _{c_t}\ga_{c_{t-1}+1}(\ldots (
 \ga_{c_1+1}(F))\ldots )]}.$$
 Hence the following useful isomorphism holds
$${\cal N}_{c_1,\ldots ,c_t}M(G)\cong {\cal N}_{c_t}M(\ldots ({\cal
N}_{c_1}M(G))\ldots ).$$ Clearly ${\cal N}_{c_1}M(G)\cong
\ga_{c_1+1}(F)/\ga_{n+c_1+1}(F)$ if the free abelian group of
rank $\sum_{i=c_1+1}^{c_1+n}\chi_i(m)$. We need the following
theorem which can be easily proved similar to the proof of the
main result in a joint paper of the first author [10].\\
\textbf{Theorem 2.10}. \textit {Let $G={\bf Z}^{(m)}\oplus {\bf
Z}_{n_1}\oplus {\bf Z}_{n_2}\oplus \ldots \oplus {\bf Z}_{n_k}$ be a
finitely generated abelian group, where $n_{i+1}|n_i$ for all $1\leq
i\leq k-1$. Then, for all $c\geq 1$
$$ {\cal N}_cM(G)\cong {\bf Z}^{(b_m)}\oplus{\bf Z}_{n_1}^{(b_{m+1}-b_m)}\oplus {\bf
Z}_{n_3}^{(b_3-b_2)}\oplus \ldots \oplus {\bf
Z}_{n_k}^{(b_{m+k}-b_{m+k-1})} \ \ \ \ ,$$ where
$b_i=\chi_{c+1}(i)$ and $X^{(m)}$ denotes the direct sum of
$m$-copies of $X$}.\\

 Now, using induction and the above facts we can present an
 explicit formula for the polynilpotent multiplier of free
 nilpotent groups as follows which extends somehow the results of
 [10,14,15,16,17,18,20].\\
\textbf{Theorem 2.11}. \textit {Let $G=\textbf{ Z}\st{n}* \textbf{
Z}\st{n}*\ldots \st{n}*\textbf{ Z}$ be the free $n$th nilpotent
group of rank $m$ and ${\cal N}_{c_1,\ldots ,c_t}$ be the variety of
polynilpotent groups of class row $(c_1\ldots ,c_t)$. If $c_1\geq
n$, then ${\cal N}_{c_1,\ldots ,c_t}M(G)$ is the free abelian group
of rank $\chi_{c_t+1}(\ldots
(\chi_{c_2+1}(\sum_{i=c_1+1}^{c_1+n}\chi_i(m)))\ldots )$}.\\
\  \  \  \\
\textbf{ACKNOWLEDGMENT}

 This research was in part supported by a grant from \textbf{IPM} (No. 84200038).

\end{document}